
\input amstex
\magnification 1200
\documentstyle{amsppt}
\topmatter
\title Ternary quartics and 3-dimensional commutative algebras
\endtitle
\author P.Katsylo, \ D.Mikhailov
\endauthor
\address Independent University of Moscow
\endaddress
\date May 5, 1994
\enddate
\thanks Research supported by Max-Planck-Institute fur Mathematik and
 Grant N MQZ000 of the International Science Foundation.
\endthanks
\abstract We find the connection between 3-dimensional
 commutative algebras with a trivial trace and plane
 quartics and its bitangents.
\endabstract
\endtopmatter
\document
\define\C{\Bbb C}
\define\p{\Bbb P}
\par
\S\bf 0. \rm A structure of a commutative algebra on $\C^3$ is called
3-dimensional algebra. Let $\Cal A$ be the set of 3-dimensional algebras.
Consider
$\Cal A$ as a linear space. Let $\Cal A_0 \subset \Cal A$ be the linear
subspace of algebras with a trivial trace. By definition,
$\eta \in \Cal A_0$ iff the contraction of the structural constants
of $\eta$ is equal to zero.
\par
By $PV$ denote the projectivization of a vector space $V$.
For $v \in V,\ \overline{v} \neq 0$ denote by $\overline{v}$
the corresponding point of the proective space $PV$.
\par
Let $\eta \in \Cal A_0$ be an algebra with a trivial trace.
Recall that an element $a \in \C^3$ is called an idempotent,
iff $a \neq 0,\ a^2 = a$. We say that an element
$\overline{a} \in P\C^3 = \p^2$ is called a generalized idempotent,
iff $a^2 = \lambda a$, there $\lambda \in \C$.
Every idempotent defines a generalized idempotent. Every generalized
idempotent $\overline{a} \in \p^2$ such that $a^2 \neq 0$ defines
uniquely the idempotent $a' \in \C^3$ such that $\overline{a} =
\overline{a'}$. Define the subscheme $X(\eta) \subset \p^2$
of the generalized idempotents by the following equation:
$$a^2 \wedge a = 0 \tag0.1$$
\par
Consider the open $SL_3$-invariant subset
$$\Cal A'_0 = \{\eta \in \Cal A_0 \mid \dim X(\eta) = 0\} \subset \Cal A_0.$$

\proclaim{Lemma 0.1} $\Cal A'_0$ is nonempty
\endproclaim

\par
Consider $\eta \in \Cal A_0$. The algebra $\eta$ defines
the quadratic mapping
$$\C^3 \longrightarrow \C^3, \ \ a \mapsto a^2.$$
This quadratic mapping defines the section $\tilde \eta$
of the vector bundle $T_{\p^2}(1)$. The scheme of zeros of
the section $\tilde \eta$ is $X(\eta)$. We have
$$\deg X(\eta) = c_2(T_{\p^2}(1)) = 7,$$
for $\eta \in \Cal A'_0$.
\par
Consider the open $SL_3$-invariant subset
$$\split
      \Cal A''_0 = & \{\eta \in \Cal A'_0 \mid X(\eta) =
      \{\overline{a_1}, \dots ,\overline{a_7}\},
      \ a^2_i \neq 0,\ 1 \le i \le 7, \\
      & every \ 3 \ points \ of \ X(\eta) \ does \ not \ lie \ on
      \ a \ line, \\
      & every \ 6 \ points \ of \ X(\eta) \ does \ not \ lie \ on
      \ a \ quadric \}
      \subset \Cal A'_0.
  \endsplit
$$
\proclaim {Lemma 0.2} $\Cal A''_0$ is nonempty.
\endproclaim
\par
Consider the rational $SL_3$-morphism
$$\phi : P\Cal A_0 \longrightarrow (\p^2)^{(7)},
  \ \eta \mapsto X(\eta).$$
\proclaim {Proposition 0.3} $\phi$ is a birational isomorphism.
\endproclaim
\par
In other words, 3-dimensional algebra in general position with
a trivial trace is uniquely (up to a scalar factor) defined
by its generalized idempotents.
\par
Fix $\eta \in \Cal A''_0$. Let $a_1, \dots , a_7$ be the idempotents
of the algebra $\eta$. Let
$$\pi = \pi(\eta) : Z = Z(\eta) \longrightarrow \p^2$$
be the blowing up of $X(\eta)$ in $\p^2$. It is well known that $Z$ is
Del Pezzo surface of degree 2. Let
$$\beta = \beta(\eta) : Z \longrightarrow \p^{2\ast}$$
be the canonical double covering with the nonsingular quartic
$Y = Y(\eta) \subset \p^{2\ast}$ as the branch locus.
\par
$SL_3$-module $S^2\Cal A_0$ contains with multiplicity 1 the isomorphic
to $S^4\C^3$ submodule. Therefore, there exists a unique (up to a scalar
factor) nontrivial quadratic $SL_3$-mapping
$$\epsilon : \Cal A_0 \longrightarrow S^4\C^3.$$
\proclaim {Lemma 0.4} $\epsilon(\eta) \neq 0$.
\endproclaim
Consider the quaternary form $\epsilon(\eta)$ on the space $\C^{3\ast}$.
This quaternary form defines quartic $Y' = Y'(\eta) \subset \p^{2\ast}$.
Consider 28 linear forms $a_1, \dots , a_7,\ (a_i - a_j)^2,\
1 \leq i < j \leq 7$ on the space $\C^{3\ast}$. These linear forms
define 28 lines $A_1, \dots , A_7,\ A_{ij} \in \p^{2\ast}$.
\proclaim {Theorem 0.5} 1) $Y$ is isomorphic to $Y'$.
   \item{2)} 28 bitangents to $Y'$ are $A_1, \dots , A_7,
            \ A_{ij},\ 1 \leq i < j \leq 7$.
\endproclaim
The authors thank B.Broer for useful discussion.
\par
\S\bf 1. \rm Let $e_1, e_2, e_3$ be the standart basis in $\C^3$, and
$x_1, x_2, x_3$ be the dual basis in $\C^{3\ast}$. The group $SL_3$
acts canonically in the space $S^a\C^3 \otimes S^b\C^{3\ast},\
a,b \geq 0$. For $a,b \geq 1$ consider the linear $SL_3$-mapping
$$\Delta = \sum {\partial \over {\partial e_i}} \otimes
  {\partial \over {\partial x_i}} : S^a\C^3 \otimes S^b\C^{3\ast}
  \longrightarrow S^{a-1}\C^3 \otimes S^{b-1}\C^{3\ast}.$$
It is well known that a representation of the group $SL_3$ in the space
$V(a,b) = \ker \delta$ is irreducible. Assume that
$V(a,0) = S^a\C^3,\ V(0,b) = S^b\C^{3\ast}$.
\par
A structure of a commutative algebra on $\C^3$ is a symmetric
bilinear mapping
$$\C^3 \times \C^3 \longrightarrow \C^3.$$
The set of such symmetric bilinear mappings is
$\C^3 \otimes S^2\C^{3\ast}$. Therefore, the linear space $\Cal A$
of 3-dimensional algebras is $\C^3 \otimes S^2\C^{3\ast}$. The contraction
of structural constants of algebras is the mapping
$$\Delta : \C^3 \otimes S^2\C^{3\ast} \longrightarrow \C^{3\ast}.$$
Therefore, the linear space $\Cal A_0$ of 3-dimensional algebras
with a trivial trace is $V(1,2)$.
\par
The decomposition
$$V(1,2) \otimes \C^3 \otimes \C^3 \simeq \C^3 \oplus 2V(0,2)
  \oplus 2V(2,1) \oplus 2V(1,3) \oplus V(3,2)$$
holds. Thus there exists a unique (up to a scalar factor)
nontrivial trilinear $SL_3$-mapping
$$\mu : V(1,2) \times \C^3 \times \C^3 \longrightarrow \C^3.$$
Let us give an explicite form of $\mu$:
$$\mu(e_{i_1} \otimes x_{j_1}x_{j_2}, e_{i_2}e_{i_3}) =
  \Delta^2(e_{i_1}e_{i_2}e_{i_3} \otimes x_{j_1}x_{j_2}).$$
The corresponding to $\eta \in V(1,2)$ algebraic structure is
the bilinear symmetric mapping
$$\mu(\eta, \cdot, \cdot) : \C^3 \times \C^3 \longrightarrow \C^3.$$
\par
\demo{Example 1.1} Consider the algebra
   $$\eta_0 = \frac{1}{4} (e_1 \otimes x_3^2 + e_2 \otimes x_1^2 +
     e_3 \otimes x_2^2).$$
   The multiplication table of $\eta_0$ is as following:
   $$e_1 \ast e_2 = e_2 \ast e_3 = e_3 \ast e_1 =0,
     \ e_1^2 = e_2,\ e_2^2 = e_3,\ e_3^2 = e_1.$$
   It can be easily checked that the subscheme $X(\eta_0) \subset \p^2$
   of the generalized idempotents of $\eta_0$ is
   $$X(\eta_0) = \{\overline{\theta e_1 + \theta^2 e_2 + \theta^4 e_3}
     \mid \theta \in \mu_7\}.$$
\enddemo
\demo {Proof of lemmas 0.1 and 0.2}
   We have $\eta_0 \in \Cal A'_0$. It follows that $\Cal A'_0$
   is nonempty. It can easily be checked that $\eta_0 \in \Cal A''_0$.
   Therefore, $\Cal A''_0$ is nonempty.
\enddemo
\demo {Proof of proposition 0.3}
   It can easily be checked that $\phi^{-1}(X(\eta_0)) =
   \overline{\eta_0}$. It follows from (0.1) that a closure of
   a fiber of $\phi$ is a linear subspace in $PV(1,2)$ of dimension
   $\leq \dim{PV(1,2)} - \dim(\p^2)^{(7)} = 0$. Hence $\phi$ is
   a birational isomorphism.
\enddemo
\S\bf 2. \rm Fix $\eta \in \Cal A''_0$. Let $a_1, \dots, a_7$
be the idempotents of $\eta$.
\par
Consider the cubic mapping
$$\psi = \psi(\eta) : \C^3 \longrightarrow \wedge^2\C^3
  \simeq \C^{3\ast}, \ a \mapsto a^2 \wedge a.$$
\proclaim {Lemma 2.1} Consider $a_i$ as a linear form
   on $\C^{3\ast}$. Let $Q_i = Q_i(\eta)$ be the corresponding
   to the cubic form $\psi^\ast(a_i)$ cubic. Then the cubic
   $Q_i$ contains $\overline{a_1}, \dots ,
   \overline{a_7}$. Moreover, $Q_i$
   contains $\overline{a_i}$ with multiplicity $\geq$ 2.
\endproclaim
\demo {Proof} It is obvious that $Q_i \ni \overline{a_1},
   \dots ,\overline{a_7}$.
   Let us prove that $Q_i$ contains $\overline{a_i}$
   with multiplicity $\geq$ 2. We have
   $$\psi^\ast(a_i) : a \mapsto a^2 \wedge a \wedge a_i
     \in \wedge^3\C^3 \simeq \C,$$
   $$\split
     \psi^\ast(a_i)(a_i + tb)
     &= (a_i + tb)^2 \wedge (a_i + tb)
     \wedge a_i = a_i^2 \wedge a_i \wedge a_i + \\
     &t((2a_i \ast b) \wedge a_i
     \wedge a_i + a_i^2 \wedge b \wedge a_i) + \dots =
     0 +t \cdot 0 + \dots
     \endsplit$$
   for any $b \in \C^3$.
\enddemo
\proclaim {Lemma 2.2} Consider $(a_i - a_j)^2, i < j$ as a linear
   form on $\C^{3\ast}$. Let $Q_{ij} = Q_{ij}(\eta)$ be the corresponding
   to a cubic form $\psi^\ast((a_i - a_j)^2)$ cubic. Then $Q_{ij}$
   is the union of the line $\langle\overline{a_i},\overline{a_j}\rangle$
   and the containing the points
   $\overline{a_1}, \dots ,\widehat{\overline{a_i}}, \dots ,
   \widehat{\overline{a_j}}, \dots , \overline{a_7}$
   quadric.
\endproclaim
\demo {Proof} It is obvious that
   $Q_{ij} \ni \overline{a_1}, \dots ,\overline{a_7}$. We have to prove
   that $Q_{ij}$ contains the line
   $\langle\overline{a_i},\overline{a_j}\rangle$. We have
   $$\psi^\ast((a_i - a_j)^2) : a \mapsto a^2 \wedge a
     \wedge (a_i - a_j)^2 \in \wedge^3\C^3 \simeq \C,$$
   $$\split
       \psi^\ast((a_i & - a_j)^2)(t_ia_i + t_ja_j) = (t_ia_i +t_ja_j)^2
       \wedge (t_ia_i + t_ja_j) \wedge (a_i - 2a_i \ast a_j + a_j) = \\
       &t_i^3a_i^2 \wedge a_i \wedge (a_i - 2a_i \ast a_j + a_j) +
       t_i^2t_j(a_i^2 \wedge a_j \wedge (a_i - 2a_i \ast a_j + a_j) + \\
       &((2a_i \ast a_j) \wedge a_i \wedge (a_i - 2a_i \ast a_j + a_j)) +
       \dots = t_i^3 \cdot 0 + t_i^2t_j \cdot 0 + \dots \equiv 0.
     \endsplit $$
\enddemo
\par
Consider the rational morphism
$$\Psi = \Psi(\eta) : \p^2 \longrightarrow \p^{2\ast},
  \ \overline{a} \mapsto \overline{\psi(a)}.$$
The rational morphism $\Psi$ is not defined on $X(\eta)$. Let
$$\CD
  \p^2  @<\pi=\pi(\eta)<<  Z=Z(\eta)  @>\beta=\beta(\eta)>>  \p^{2\ast}
  \endCD $$
be the regularization of $\Psi$. It is well known that $Z$ is Del Pezzo
surface of degree 2, $\pi$ is the blowing up of the seven points
$\overline{a_1}, \dots ,\overline{a_7}$ in
$\p^2, \beta : Z \longrightarrow \p^{2\ast}$ is double covering with
the nonsingular quartic $Y = Y(\eta) \subset \p^{2\ast}$
as the branch locus.
\demo {Proof of lemma 0.4 and theorem 0.5} Consider the nontrivial
   homogeneous (of degree 6) $SL_3$-mappings
   $$\gamma_i : V(1,2) \longrightarrow V(0,12) = S^{12}\C^{3\ast},$$
   $$\gamma_1 : \eta \mapsto \psi^\ast(\epsilon(\eta)),$$
   $$\gamma_2 : \eta \mapsto
     (\det(\frac{\partial\psi^\ast(e_i)}{\partial x_j}))^2.$$
   The decomposition
   $$S^6V(1,2) \simeq \dots \oplus V(1,12) \oplus \dots$$
   holds. Thus $\gamma_1 = c\gamma_2$, there $c \neq 0$. This implies
   lemma 0.4 and statement 1) of theorem 0.5.
   \par
   Statement 2) of theorem 0.5 is a corollary of lemmas 2.1 and 2.2
   (see \cite{1})
\enddemo

\enddocument